\documentclass[11pt,letterpaper,reqno]{amsart}
\usepackage[left=30mm,top=30mm,right=30mm,bottom=30mm]{geometry}
\usepackage{epsfig}
\usepackage{amssymb}
\usepackage{amsthm}
\usepackage{subcaption}
\usepackage{multicol}
\usepackage{times}
\usepackage{color}
\usepackage{enumerate}
\usepackage{graphicx}
\usepackage{nicefrac}
\usepackage{xcolor}
\usepackage{tikz}
\usepackage[all]{xy}
\usepackage{comment}
\usepackage{hyperref}
\hypersetup{
 colorlinks=true, 
 linkcolor=blue, 
 citecolor=blue, 
 filecolor=blue, 
 urlcolor=blue 
}
\usepackage[nameinlink,capitalize,noabbrev]{cleveref}

\newtheorem{theorem}{Theorem}
\newtheorem*{theorem*}{Theorem}

\newtheorem{lemma}[theorem]{Lemma}

\newtheorem*{claim*}{Claim}

{\theoremstyle{remark}
 \newtheorem{remark}[theorem]{Remark}}
 \newtheorem*{remark*}{Remark}
{\theoremstyle{definition}

}

\newcommand{\Z}{\mathbb Z}

\newcommand{\R}{\mathbb R}

\newcommand{\K}{\mathbf{k}}

\newcommand{\hs}{\mathcal{H}(S)}

\begin{document}

\title[Nash blowup fails to resolve singularities]{Nash blowup fails to resolve singularities \\ in dimensions four and higher}

\author[{F. Castillo}]{Federico Castillo}
\address{Facultad de Matemáticas, Pontificia Universidad Católica de Chile, Santiago, Chile}  
\email{federico.castillo@uc.cl}

 \author[D. Duarte]{Daniel Duarte}
\address{Centro de Ciencias Matematicas, UNAM, Morelia, M\'exico}  \email{adduarte@matmor.unam.mx}

\author[M. Leyton-Alvarez]{Maximiliano Leyton-\'Alvarez} %
\address{Instituto de Matem\'aticas, Universidad de Talca, Talca, Chile} %
\email{mleyton@utalca.cl}

\author[A. Liendo]{Alvaro Liendo} %
\address{Instituto de Matem\'aticas, Universidad de Talca, Talca, Chile} %
\email{aliendo@utalca.cl}

\date{\today}

\thanks{
 \mbox{\hspace{11pt}}{\it Key words}: Resolution of singularities, Nash blowup, Toric varieties.\\
 \mbox{\hspace{11pt}} The first author was partially supported by Fondecyt project 1221133. The second author was partially supported by CONAHCYT project CF-2023-G-33 and PAPIIT grant IN117523. The third author was partially supported by Fondecyt project 1221535. The fourth author was partially supported by Fondecyt project 1240101.  Finally, the first, third, and fourth authors were partially supported by Fondecyt Exploraci\'on project 13250049. 
 }

\begin{abstract}
In this paper we show that iterating Nash blowups or normalized Nash blowups does not resolve the singularities of algebraic varieties of dimension four or higher over an algebraically closed field of arbitrary characteristic.
\end{abstract}

\dedicatory{Dedicated to the memory of Gerard Gonzalez-Sprinberg}

\maketitle

\section*{Introduction}

Let $X \subseteq \K^n$ be an equidimensional algebraic variety of dimension $d$, where $\K$ is an algebraically closed field. Consider the Gauss map:
\begin{align}\label{eq:gauss}
\Phi\colon X \setminus \operatorname{Sing}(X) \to \operatorname{Grass}(d, n) \quad \text{defined by} \quad x \mapsto T_x X\,,
\end{align}
where $\operatorname{Grass}(d, n)$ denotes the Grassmannian of $d$-dimensional vector spaces in $\K^n$, and $T_x X$ is the tangent space to $X$ at $x$. Let $X^*$ be the Zariski closure of the graph of $\Phi$ and $\nu\colon X^* \to X$ be the composition of the inclusion $X^*\hookrightarrow X\times\operatorname{Grass}(d, n)$ and the projection onto the first coordinate.
The morphism $\nu$ is a proper birational map that is an isomorphism over $X \setminus \operatorname{Sing}(X)$. The pair $(X^*,\nu)$ is called the \emph{Nash blowup} of $X$. Letting now $\eta\colon \overline{X^*}\to X^*$ be the normalization map, the pair $(\overline{X^*},\eta\circ\nu)$ is called the \emph{normalized Nash blowup} of $X$. 

It was proposed by J. Nash to resolve singularities by iterating the Nash blowup \cite{Sp90}. This question was also previously proposed by J. G. Semple in the 1950's \cite{Se54}. As a starting point for the study of this question, A. Nobile proved that the Nash blowup is an isomorphism if and only if the variety is smooth, over fields of characteristic zero \cite{Nob75}. This result served as a motivation and several partial results were obtained over the years giving an affirmative answer for some families of varieties over fields of characteristic zero \cite{Nob75,Reb,EGZ,GoTe14,DuarSurf,DG}. 

Possibly inspired by Zariski's method of resolution of surfaces, which consists of normalized point blowing up, G. Gonz\'alez-Sprinberg introduced a variant to Nash's question that consists in iterating the normalized Nash blowup \cite{GS1,GS0}. This variant received a lot of attention and significant results were obtained in characteristic zero by G. Gonz\'alez-Sprinberg, H. Hironaka, and M. Spivakovsky, among others \cite{GS1,GS2,Hi83,Sp90,GS3,Ataetal,GrMi12}. 

In the case of positive characteristic fields, the study of Nash blowups was discouraged by the following example of A. Nobile \cite{Nob75}: the Nash blowup of the curve $x^p - y^q = 0$ in $\K^2$ is isomorphic to the curve when the characteristic of $\K$ is $p$ or $q$, where $p$ and $q$ are distinct prime numbers. This bad behaviour was overcome by adding the condition of normality. D. Duarte and L. N\'u\~nez Betancourt showed that, for normal varieties over fields of positive characteristic, the Nash blowup is an isomorphism if and only if the variety is smooth \cite{DuNu22}. Hence, the question on the resolution properties of normalized Nash blowup was reconsidered over fields of prime characteristic as well. In this context, some new results were recently obtained \cite{DJNB,DDR,CDLLcharfree}.

In summary, there are longstanding conjectures on whether the iteration of Nash blowups resolves singularities in characteristic zero and whether normalized Nash blowup resolves singularities in arbitrary characteristics. If true, these procedures could provide a canonical algorithm for resolution of singularities. 

The aim of this paper is to disprove both the Nash blowup conjecture and the normalized Nash blowup conjecture. More precisely, we prove the following theorem.

\begin{theorem*}\label{main}
For every $d\geq 4$ and every algebraically closed field $\K$, there exists a normal singular affine algebraic variety $X$ of dimension $d$ over $\K$ such that:
\begin{enumerate}
    \item[$(i)$] If $\operatorname{char}(\K)=0$, then the Nash blowup of $X$ and the normalized Nash blowup of $X$ contain an open affine subset isomorphic to $X$.
     \item[$(ii)$] If $\operatorname{char}(\K)$ is positive and different from $3$, then the normalized Nash blowup of $X$ contains an open affine subset isomorphic to $X$.
    \item[$(iii)$] If $\operatorname{char}(\K)=3$, then the second iteration of the normalized Nash blowup of $X$ contains an open affine subset isomorphic to $X$.
\end{enumerate}
Consequently, in every dimension $d \geq 4$ and every characteristic, iterating the Nash blowup or the normalized Nash blowup fails to resolve singularities.

\end{theorem*}

To prove the theorem, we construct explicit examples of normal singular affine toric varieties corresponding to each of the described cases. In all characteristics, these examples arise as affine charts appearing in an iteration of the normalized Nash blowup of the toric variety associated to the semigroup $S = \omega \cap \mathbb{Z}^4$, where $\omega \subset \mathbb{R}^4$ is the cone generated by the columns of the matrix
$$\left[\begin{array}{rrrr}
1 & 0 & 0 & 1 \\
0 & 1 & 0 & 1 \\
0 & 0 & 1 & 1 \\
0 & 0 & 0 & 5
\end{array}\right]\, .$$

The study of Nash blowups of toric varieties over algebraically closed fields of characteristic zero was initiated by G. Gonz\'alez-Sprinberg in the seminal paper \cite{GS1}. He gave a combinatorial description of the normalized Nash blowup of a normal toric variety by using the so-called logarithmic Jacobian ideal. This result was later improved by removing the condition of normality by P. Gonz\'alez and B. Teissier \cite{GoTe14} (see also \cite{LJ-R,GrMi12}). Finally, the positive characteristic versions of these results were obtained by D. Duarte, J. Jeffries, and L. N\'{u}\~nez-Betancourt by means of an analog of the logarithmic Jacobian ideal that takes into account the characteristic \cite{DJNB}. These combinatorial descriptions for the Nash blowup or the normalized Nash blowup of a toric variety, both in characteristic zero and prime, are our main tools to compute the varieties mentioned in the theorem. 

The theorem is stated over algebraically closed fields. However, in characteristic zero there exists an alternative approach that allows us to prove it without assuming algebraic closedness. For details, see \cref{rem:non-alg-closed}.

Finally, we would like to highlight that a significant portion of the work leading to the results in this paper involved extensive computational experimentation using SageMath \cite{sagemath}. Further details of these computations are presented in \cref{sec:computer}. The computational approach to studying Nash blowups of toric varieties has also been employed in previous works, such as \cite{Ataetal}. 

In view of the theorem and the known results on Nash blowups and normalized Nash blowups, the remaining cases to completely answer the conjectures on the resolution properties of Nash blowups are:

\begin{itemize}
\item[(A)] Iteration of the Nash blowup of algebraic varieties of dimension 2 and 3 over algebraically closed fields of characteristic zero.
\item[(B)] Iteration of the normalized Nash blowup of normal algebraic varieties of dimension 3 over algebraically closed fields of characteristic zero.
\item[(C)] Iteration of the normalized Nash blowup of normal algebraic varieties of dimension 2 and 3 over algebraically closed fields of prime characteristic.
\end{itemize}

In the case of toric varieties, the previous list is the same except in problem (C), where the remaining case is normal toric varieties of dimension 3 \cite{DJNB}.

\subsection*{Acknowledgments}

This collaboration began at the \href{https://sites.google.com/view/agrega0/home}{AGREGA} workshop, which took place at Universidad de Talca in January 2024 where the second named author presented the Nash blowup conjecture. We extend our gratitude to the institution for its support and hospitality. We would also like to express our sincere gratitude to the editors and the anonymous referees for their invaluable suggestions to greatly improve our manuscript.

\section{The Nash blowup of an affine toric variety}
\label{sec:toric}

A toric variety $X$ is a variety endowed with a faithful regular action of an algebraic torus having an open orbit. Toric varieties admit a well-known combinatorial description that we recall now, for details, see \cite{fulton1993introduction,oda1983convex,Sturm,cox2011toric}. In contrast to some references, but in line with \cite{Sturm,cox2011toric}, we do not require toric varieties to be normal.

Let $S$ be an affine semigroup, i.e., a finitely generated semigroup with identity element that can be embedded in a free abelian group. In the sequel, without loss of generality, we assume that $S\subset M=\Z^d$ and that the group generated by $S$ is $M$. Let now $M_\R=M\otimes_\Z\R$ and let  $\omega\subset M_\R$ be the polyhedral cone generated by $S$, i.e.,
$$\omega = \operatorname{Cone}(S) = \left\{ \sum_{u \in F}  \lambda_u u \mid F\subset S \mbox{ finite and } \lambda_u \geq 0 \right\} \subseteq \R^d\,.$$

We say that $S$ is saturated if $S=\omega\cap M$. We say that the semigroup $S$ is pointed if $\omega$ is strongly convex, i.e., if $\{0\}$ is the only vector space contained in $\omega$. The Hilbert basis of a pointed affine semigroup $\mathcal{H}(S)$ is its minimal generating set. It consists of elements $h\in S$ that cannot be written as a sum $h=m+m'$ with $m,m'\in S$ and $m,m'\neq 0$.

Given an affine semigroup $S$, we define the semigroup algebra
\[
\K[S]=\bigoplus_{u\in S}\K\cdot\chi^{u},\quad\mbox{with}\quad\chi^{0}=1,\mbox{ and }\chi^{u}\cdot\chi^{u'}=\chi^{u+u'}, \mbox{ for all } u,u'\in S\,.
\]
The affine variety $X(S):=\operatorname{Spec}\K[S]$ is a toric variety. Moreover, the toric variety $X(S)$  is  normal if and only if $S$ is saturated \cite[Theorem~1.3.5]{cox2011toric}.

\medskip

As mentioned in the introduction, the Nash blowup and the normalized Nash blowup of a toric variety over fields of arbitrary characteristic can be described combinatorially. We now recall this description following the account given in \cite[Section~1.9.2]{Sp20}.

Let $X(S)$ be the affine toric variety given by the pointed semigroup $S\subset \Z^d$ and let $\mathcal{H}(S)=\{h_1,\dots,h_r\}$ be its Hilbert basis. For a collection of $d$ elements $\{h_{i_1},\dots,h_{i_d}\}\subset\hs$, we define the matrix $(h_{i_1} \cdots h_{i_d})$ whose columns are the vectors $h_{i_j}$. Letting $p\geq0$ be the characteristic of the base field $\K$ we denote 
$$\operatorname{det}_p(h_{i_1}\cdots h_{i_d})=
\begin{cases}
  \det(h_{i_1}\cdots h_{i_d}) & \mbox{if } p=0 \\
  \det(h_{i_1}\cdots h_{i_d}) \mod p & \mbox{if } p>0
\end{cases}\,.$$

The affine charts of the Nash blowup of $X$ are indexed by the subsets $A=\{h_{i_1},\dots,h_{i_d}\}$ of $\hs$ such that $$\operatorname{det}_p(h_{i_1}\cdots h_{i_d})\neq 0\,.$$

Let \(A\) be such a subset of \(\mathcal{H}\), and let \(h \in A\).
Without loss of generality, up to reordering the indices, we may and will assume $A=\{h_1,\dots,h_d\}$ and $h=h_1$. We let

\begin{align} \label{eq:Ga}
\mathcal{G}_A(h)=\Big\{g-h\mid g\in \hs\setminus A \mbox{ and } \operatorname{det}_p(g\,h_2\cdots h_d)\neq 0 \Big\}.
\end{align}
Finally, letting
$\mathcal{G}_A=\hs\cup \mathcal{G}_A(h_1)\cup\dots\cup \mathcal{G}_A(h_d)$, we let  $S_A$ be the semigroup generated in $M$ by $\mathcal{G}_A$. Now, the collection of affine toric varieties $X(S_A)$, for all $A=\{h_{i_1},\dots,h_{i_d}\}$ with 
$$\operatorname{det}_p(h_{i_1}\cdots h_{i_d})\neq 0 \quad \mbox{and}\quad S_A\mbox{ pointed}, $$
provides a set of covering affine charts of the Nash blowup of $X(S)$ (see \cite[Propositions 32 and 60]{GoTe14} for the characteristic zero case and \cite[Proposition 32]{GoTe14} and \cite[Theorem 1.9]{DJNB} for the prime characteristic case).

Let us now recall the combinatorial description of the normalized Nash blowup of a toric variety, which is a consequence of the previous description.
Consider the following polyhedron associated to the toric variety $X(S)$:
\begin{align} \label{eq:norm-Nash}
\mathcal{N}_p(S)=\operatorname{Conv}\left(\{(h_{i_1}+\dots+h_{i_d})+\omega\mid h_{i_j}\in \mathcal{H}(S)\mbox{ and }\operatorname{det}_p(h_{i_1}\cdots h_{i_d})\neq 0\}\right).
\end{align}
Here $\operatorname{Conv}(\cdot)$ denotes the convex hull of a set. Then a set of covering affine charts of the normalized Nash blowup of $X(S)$ is given by $X(S_v)$, where $v$ runs over all the vertices of the polyhedron $\mathcal{N}_p(S)$, $\omega_v = \operatorname{Cone}\left(\mathcal{N}_p(S)-v\right)\subset M_\R$, and $S_v=\omega_v\cap M$.

Let us briefly recall the relation among the semigroups $S_A$ and $S_v$ coming from the Nash blowup and the normalized Nash blowup, respectively. Let $A=\{h_{i_1},\dots,h_{i_d}\}\subset\hs$ be such that $\det_p(h_{i_1}\cdots h_{i_d})\neq0$. Denote as $\overline{S_A}$ the saturation of $S_A$, that is, $\overline{S_A}=\operatorname{Cone}(S_A)\cap M$. Now let $v=h_{i_1}+\dots+h_{i_d}$. Then $S_A$ being pointed is equivalent to $v$ being a vertex of $\mathcal{N}_p(S)$. Since it is known that normalizing a toric variety corresponds to saturating the semigroup we obtain
$$\operatorname{Cone}\left(\mathcal{N}_p(S)-v\right)=\operatorname{Cone}\left(\mathcal{G}_A\right) \quad\mbox{ and }\quad S_v=\overline{S_A}\,.$$

In the next section we will use these descriptions and remarks to provide the examples stated in the theorem.

\section{Proof of the Theorem} \label{sec:proof-counter}

We first state the following lemma that follows directly from the Gauss map in \eqref{eq:gauss}.
\begin{lemma} \label{lem:product}
    Let $X$ be an affine variety and let $X^*\to X$ be its Nash blowup. Then  $X^*\times\K^r\to X\times \K^r$ is the Nash blowup of $X\times \K^r$.
\end{lemma}

By \cref{lem:product}, it is enough to show the theorem for dimension $d=4$. In the sequel, we restrict ourselves to this case. As stated in the introduction, to prove the theorem we will exhibit normal singular affine toric varieties $X$ fulfilling the conditions of the theorem. Let us first prove the case where the characteristic of the base field $\K$ is 0.

Let $M=\Z^4$ and consider the cone $\omega$ in $M_\R=\R^4$ generated by the columns of the matrix
\begin{align} \label{eq:count-char-0}
B=\left[\begin{array}{rrrrrr}
1 & 0 & 0 & 0 & 2 & 1 \\
0 & 1 & 0 & 0 & 3 & 3 \\
0 & 0 & 1 & 0 & -2 & -1 \\
0 & 0 & 0 & 1 & -1 & -1
\end{array}\right]\,.
\end{align}
Let $X(S)$ be the toric variety defined by $S = \omega \cap M$. Each column of the matrix $B$ corresponds to the primitive generator of a ray of the cone $\omega$. We denote by $h_i$ the vector corresponding to the $i$-th column of $B$, with $i \in \{1, \dots, 6\}$.

\begin{lemma}\label{lem S char 0}
The Hilbert basis of $S$ is $\mathcal{H}(S)=\{h_1,h_2,\dots,h_6,h_7\}$, where $h_7=(1, 2, -1, 0)$.
\end{lemma}

\begin{figure}[ht]
\begin{center}
\begin{tikzpicture}[scale=0.75]
\newcommand{\scl}{0.7}
\newcommand{\xiooi}{0*\scl}	\newcommand{\yiooi}{6*\scl}
\newcommand{\xooii}{3*\scl}	\newcommand{\yooii}{4*\scl}
\newcommand{\xioio}{-1*\scl}	\newcommand{\yioio}{4*\scl}
\newcommand{\xoioi}{1*\scl}	\newcommand{\yoioi}{2*\scl}
\newcommand{\xiioo}{-3*\scl}	\newcommand{\yiioo}{2*\scl}
\newcommand{\xoiio}{0*\scl}	\newcommand{\yoiio}{0*\scl}
\coordinate (oiio) at (\xoiio,\yoiio);
\coordinate (ioio) at (\xioio,\yioio);
\coordinate (ooii) at (\xooii,\yooii);
\coordinate (iioo) at (\xiioo,\yiioo);
\coordinate (oioi) at (\xoioi,\yoioi);
\coordinate (iooi) at (\xiooi,\yiooi);
\foreach \coo in {(oiio),(iooi),(iioo),(ooii)} \draw[dotted] (ioio)--\coo;

\draw[dashed] (oiio)--(iooi);

\draw[thick] (oiio)--(iioo)--(iooi)--(ooii)--(oioi)--(iioo);
\draw[thick] (oioi)--(oiio)--(ooii) (oioi)--(iooi);

\draw[fill=black] (\xooii,\yooii) circle(.1); \node at (\xooii+.75,\yioio) {\footnotesize $h_6$};
\draw[fill=black] (\xoioi,\yoioi) circle(.1); \node at (\xoioi+.75,\yoioi) {\footnotesize $h_5$};
\draw[fill=black] (\xoiio,\yoiio) circle(.1); \node at (\xoiio,\yoiio-.5) {\footnotesize $h_2$};
\draw[fill=black] (\xioio,\yioio) circle(.1); \node at (\xioio-.75,\yioio) {\footnotesize $h_3$};
\draw[fill=black] (\xiioo,\yiioo) circle(.1); \node at (\xiioo-.75,\yiioo) {\footnotesize $h_4$};
\draw[fill=black] (\xiooi,\yiooi) circle(.1); \node at (\xiooi,\yiooi+.5) {\footnotesize $h_1$};

\draw[fill=red] (-.5,.75) circle(.1); \node[above] at (-.5,.75) {\footnotesize \textcolor{red}{$h_7$}};

\end{tikzpicture}
\end{center}
\caption{The combinatorics of an affine slice of $\omega$. The vector $h_7$ lies on the face $\operatorname{Cone}(h_2,h_4,h_5)$.}
\label{fig:octahedron}
\end{figure}
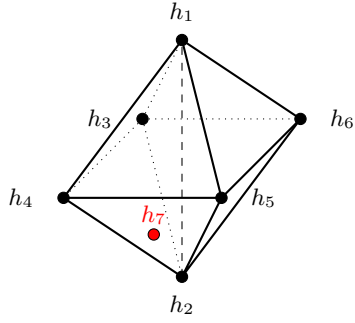

\begin{proof}
We first compute a simplicial subdivision of $\omega$. Observe that $\omega$ is a cone over a three dimensional octahedron. By choosing a pair of non-adjacent vertices, an octahedron can be triangulated with four simplices all containing the chosen pair. In the case of $\omega$, we have that $\operatorname{Cone}(h_1,h_2)$ is not a face of the cone, hence we can use it to create a simplicial subdivision of $\omega$. The following are the resulting four cones (see Figure \ref{fig:octahedron}).
\begin{align*}
\sigma_1=\operatorname{Cone}(h_1,h_2,h_3,h_4), \quad
\sigma_2&=\operatorname{Cone}(h_1,h_2,h_3,h_6), \\ \sigma_3=\operatorname{Cone}(h_1,h_2,h_5,h_6), \quad
\sigma_4&=\operatorname{Cone}(h_1,h_2,h_4,h_5)\,.  
\end{align*}
By \cite[Section VII, Theorem 2.5]{barvinok2002course} the determinant of a rational simplicial cone $\sigma$ is the same as the number of integral elements in its fundamental parallelogram $\Pi$, and the ray generators of the cone together with the elements in $\Pi$ generate the semigroup $\sigma\cap M$.

By computing determinants one can check that each of $\sigma_1, \sigma_2, $ and $\sigma_3$ is regular. The determinant of $\sigma_4$ is two, which means that there exists a unique non-zero integer vector inside its fundamental parallelogram.
This point is $h_7 = \frac12 h_2 + \frac12 h_4 + \frac12 h_5 \in \operatorname{Cone}(h_2,h_4,h_5)$ which is a face of $\omega$.
We conclude that the set $\mathcal{H}(S)$ generates the semigroup $S$.
It remains to show that it is minimal.
Every primitive ray generator is in the Hilbert basis so $\{h_1,h_2,h_3,h_4,h_5,h_6\} \subset \mathcal{H}(S)$. Since the element $h_7$ is contained in the interior of a simplicial three-dimensional face, its unique representation as a positive combination of elements in $\{h_1,h_2,h_3,h_4,h_5,h_6\}$ is   $h_7 = \frac12 h_2 + \frac12 h_4 + \frac12 h_5$ and thus $h_7$ is part of the Hilbert basis.
\end{proof}

With this lemma at hand, we can continue building the example that proves the theorem in characteristic zero.

Notice that $\operatorname{det}_0(h_1\,h_2\,h_3\,h_5)=-1$. To prove the theorem, we will show that the covering chart of the Nash blowup of $X(S)$ corresponding to $A=\{h_1,h_2,h_3,h_5\}\subset \mathcal{H}(S)$ is isomorphic to $X(S)$. 

We have $\hs\setminus A=\{h_4,h_6,h_7\}$. There are 12 relevant determinants to be computed that we exhibit below (see (\ref{eq:Ga})). At the beginning of each line we show the element of $A$ that is omitted. 
\begin{align} \label{eq:dets}
\begin{array}{llll}
h_1:&\operatorname{det}_0(h_4\,h_2\,h_3\,h_5)=-2, & \operatorname{det}_0(h_6\,h_2\,h_3\,h_5)=1, &
\operatorname{det}_0(h_7\,h_2\,h_3\,h_5)=-1, 
\\
h_2:&\operatorname{det}_0(h_4\,h_1\,h_3\,h_5)=3, & \operatorname{det}_0(h_6\,h_1\,h_3\,h_5)=0, &
\operatorname{det}_0(h_7\,h_1\,h_3\,h_5)=2,
\\
h_3:&\operatorname{det}_0(h_4\,h_1\,h_2\,h_5)=2, & \operatorname{det}_0(h_6\,h_1\,h_2\,h_5)=-1, &
\operatorname{det}_0(h_7\,h_1\,h_2\,h_5)=1,
\\
h_5:&\operatorname{det}_0(h_4\,h_1\,h_2\,h_3)=-1, & \operatorname{det}_0(h_6\,h_1\,h_2\,h_3)=1, &
\operatorname{det}_0(h_7\,h_1\,h_2\,h_3)=0.
\end{array}
\end{align}

We conclude that $\mathcal{G}_A$ is
\begin{align*}
\mathcal{G}_A=\hs&\cup\{h_4-h_1,h_6-h_1,h_7-h_1\}\cup\{h_4-h_2,h_7-h_2\}\\ &\cup\{h_4-h_3,h_6-h_3,h_7-h_3\}\cup\{h_4-h_5,h_6-h_5\}\,.
\end{align*}
Recall that $S_A$ denotes the semigroup generated by $\mathcal{G}_A$. The following subset $H$ of $\mathcal{G}_A$ is a generating set for $S_A$ (in fact, it is a Hilbert basis of $S_A$, as we will see later in the proof)
$$H=\Big\{h_1,h_4-h_2,h_7-h_2,h_4-h_3,h_6-h_3,h_4-h_5,h_6-h_5\Big\}\,.$$
Indeed, for all the other elements in $\mathcal{G}_A$ we have
\begin{align*}
h_4&=(h_4-h_3)+h_3, & h_2&=(h_7-h_2)+(h_6-h_5), \\
h_6&=(h_6-h_3)+h_3, & h_7-h_3&=(h_6-h_3)+(h_4-h_2), \\
h_7&=(h_7-h_2)+h_2, & h_4-h_1&=(h_6-h_5)+(h_4-h_3),\\
h_3&=(h_6-h_5)+h_1,& h_6-h_1&=(h_6-h_3)+(h_6-h_5),\\
h_5&=(h_6-h_3)+h_1,& h_7-h_1&=(h_6-h_3)+(h_6-h_5)+(h_4-h_2)\,.
\end{align*}
To conclude, let $U\colon M\to M$ be the automorphism of $M=\Z^4$ given by the unimodular matrix
$$U=\left[\begin{array}{rrrr}
-1 & 0 & -2 & 1 \\
0 & -1 & -3 & 0 \\
1 & 0 & 2 & 0 \\
0 & 1 & 2 & 0
\end{array}\right]\,.$$
A straightforward verification yields that $U$ induces a bijection from $\mathcal{H}(S)$ to $H$. Indeed,
\begin{align*}
U(h_1)=h_6-h_5, \qquad&
U(h_2)=h_4-h_2, \\
U(h_3)=h_4-h_5, \qquad&
U(h_4)=h_1, \\
U(h_5)=h_6-h_3, \qquad&
U(h_6)=h_4-h_3, \\
U(h_7)=h_7-h_2\,.\qquad &
\end{align*}
This yields that $S_A$ is isomorphic to $S$. Since $S$ is pointed, the same holds for $S_A$ and so $X(S_A)\simeq X(S)$ is an affine chart of the Nash blowup of $X(S)$. Furthermore, since we started with a normal variety $X(S)$ and so the isomorphic chart $X(S_A)$ is also normal, this example also shows the theorem for normalized Nash blowup in characteristic zero. This achieves the proof of the theorem in the case where the characteristic of the base field $\K$ is $0$. 

\begin{remark}
Let $\K[x_1,\dots,x_7]$ be the polynomial ring in $7$ variables. The map  $\K[x_1,\dots,x_7]\to \K[S]$ given by $x_i\mapsto \chi^{h_i}$ is surjective. Its kernel $I$  is generated by 
\begin{center}
$\begin{array}{lll}
x_1x_6-x_3x_5, &
x_4x_6-x_2x_7, &
x_7^2-x_2x_4x_5,\\
x_6x_7-x_2^2x_5, &
x_3x_7-x_1x_2^2, &
x_3x_4x_5-x_1x_2x_7.
\end{array}$
\end{center}
Hence, the toric variety $X(S)$ whose Nash blowup contains an affine chart isomorphic to itself is embedded in $\K^7$ as the zero-set of the ideal $I$. The singular locus of $X(S)$ in coordinates $x_1, \dots, x_7$ is the union of the following two-dimensional linear spaces: $\{x_2=x_4=x_5=x_6=x_7=0\}$, $\{x_2=x_3=x_5=x_6=x_7=0\}$, and $\{x_2=x_3=x_4=x_6=x_7=0\}$.
\end{remark}

\subsection*{Characteristic $p\neq 2,3$}

In the case where the characteristic is different from $2$ and $3$ the same variety $X=X(S)$ proves the theorem. Indeed, by \eqref{eq:dets} we have $$\operatorname{det}_p(h_{i_1}\,h_{i_2}\,h_{i_3}\,h_{i_4})= 0 \quad\mbox{if and only if}\quad \operatorname{det}_0(h_{i_1}\,h_{i_2}\,h_{i_3}\,h_{i_4})=0\,.$$
Hence, the starting elements for the description of the Nash blowup of $X$ in characteristic $p\neq2,3$ are exactly the same as those for characteristic zero. Thus the resulting semigroups giving the affine charts of the normalized Nash blowup of $X$ are also the same.

For the remaining characteristics $p=2$ and $p=3$, the iteration of the normalized Nash blowup does indeed resolve the singularities of $X$ so we need to provide other varieties to deal with these two cases.

\subsection*{Characteristic $p=2$}

We provide now an example of an affine variety in characteristic $2$ proving the theorem. The strategy of proof is the same as in the case of characteristic 0, so we omit the details and show only the key steps.

Let $\K$ be an algebraically closed field of characteristic $2$ and let $M=\Z^4$. Consider the cone $\omega_2$ in $M_\R=\R^4$ generated by the columns of the matrix
\begin{align} \label{eq:count-char-2}
B_2=\left[\begin{array}{rrrrr}
1 & 0 & 1 & 0 & 0 \\
0 & 1 & 1 & 0 & 2 \\
0 & 0 & 2 & 0 & 2 \\
0 & 0 & 0 & 1 & 1
\end{array}\right]\,.
\end{align}
We define also $X = X(S)$ where $S = \omega_2 \cap M$. Each column of the matrix $B_2$ corresponds to the primitive generator of a ray of the cone $\omega_2$. We denote by $h_i$ the vector corresponding to the $i$-th column of $B_2$, with $i \in \{1, \dots, 5\}$.

The Hilbert basis of $S$ is $\hs=\{h_1,\dots,h_5,h_6,h_7\}$, where $h_6=(1, 1, 1, 0)$ and $h_7=(0, 1, 1, 1)$. Notice that $\operatorname{det}_2(h_1\,h_2\,h_4\,h_7)=1$. We will show that the affine chart $X(\overline{S_A})$ of the normalized Nash blowup corresponding to the subset $A=\{h_1,h_2,h_4,h_7\}\subset \hs$ is isomorphic to $X$.
The set $\mathcal{G}_A$ is given by (recall that the determinants are taken modulo $p=2$)
\begin{align*}
\mathcal{G}_A=\hs\cup\{h_3-h_1,h_6-h_1\}\cup\{h_3-h_2\}\cup\{h_5-h_4,h_6-h_4\}\cup\{h_6-h_7\}\,.
\end{align*}
Let $H$ be the set  
$$H=\Big\{h_2,h_4,h_3-h_2,h_3-h_1,h_6-h_1,h_6-h_7,(1,0,1,0)\Big\}\,.$$
The set $H$ generates a semigroup $R$ such that $S_A\subseteq R\subseteq \overline{S_A}$. Indeed, a straightforward verification yields that every element of $\mathcal{G}_A$ is contained in $R$. Moreover, every element of $H$ is contained in $\mathcal{G}_A$ except $(1,0,1,0)$, which satisfies 
$$2(1,0,1,0)=(1,0,2,0)+(1,0,0,0)=(h_3-h_2)+h_1\in S_A\,.$$
So $(1,0,1,0)\in \overline{S_A}$. In particular, $\overline{R}=\overline{S_A}$. To conclude, let $U\colon M\to M$ be the automorphism of $M=\Z^4$ given by the unimodular matrix
$$U=\left[\begin{array}{rrrr}
1 & 0 & 0 & 0 \\
0 & 0 & 0 & 1 \\
2 & 0 & -1 & 2 \\
0 & 1 & -1 & 0
\end{array}\right]\,.$$
A straightforward verification yields that $U$ induces a bijection from $\mathcal{H}(S)$ to $H$. In particular, $S\cong R$. Since $S$ is saturated, the same holds for $R$. This yields that $S\cong R=\overline{R}=\overline{S_A}$. Since $S$ is pointed, the same holds for $\overline{S_A}$ and so $X(\overline{S_A})\simeq X$ is an affine chart of the normalized Nash blowup of $X$, concluding the proof of the theorem in the case where the characteristic of the base field $\K$ is $p=2$.

\subsection*{Characteristic $p=3$}

We provide now an example of an affine toric variety in characteristic $3$ proving the theorem. In this case, we could not find a variety whose normalized Nash blowup has a chart isomorphic to the variety. But we provide an example where this situation is obtained after 2 iterations of the normalized Nash blowup.

Let $\K$ be an algebraically closed field of characteristic $3$ and let $M=\Z^4$. Consider the cone $\omega_{3,0}$ in $M_\R=\R^4$ generated by the columns of the matrix
\begin{align} \label{eq:count-char-3}
B_{3,0}=\left[\begin{array}{rrrrr}
1 & 0 & 1 & 0 & 0 \\
0 & 1 & 2 & 0 & 3 \\
0 & 0 & 3 & 0 & 3 \\
0 & 0 & 0 & 1 & 1
\end{array}\right]\,.
\end{align}
We define $X = X(S)$ where $S = \omega_{3,0} \cap M$. The Hilbert basis of $S$ is given by $\mathcal{H}(S)=\{h_1,\dots,h_9\}$, where $h_i$ corresponds to the $i$-th column of the following matrix
$$\left[\begin{array}{rrrrrrrrr}
1 & 0 & 1 & 0 & 0 & 0 & 1 & 1 & 0 \\
0 & 1 & 2 & 0 & 3 & 2 & 1 & 2 & 1 \\
0 & 0 & 3 & 0 & 3 & 2 & 1 & 2 & 1 \\
0 & 0 & 0 & 1 & 1 & 1 & 0 & 0 & 1
\end{array}\right].
$$
Notice that $\operatorname{det}_3(h_1\,h_2\,h_4\,h_9)=2$. Let $A=\{h_1,h_2,h_4,h_9\}$. 
The semigroup $S_A$ generated by $\mathcal{G}_A$ is pointed and so $X_1 = X(\overline{S_A})$ is an affine chart of the normalized Nash blowup of $X$. We denote by $\omega_{3,1}$ the cone generated by $\overline{S_A}$ in $M_\R$. Then, the cone $\omega_{3,1}$ is generated by the columns of the matrix
$$B_{3,1}=\left[\begin{array}{rrrrr}
0 & 0 & 0 & 1 & 1 \\
0 & 1 & 2 & 0 & 1 \\
0 & 0 & 3 & 0 & 3 \\
1 & 0 & 0 & -1 & 0
\end{array}\right]\,.$$
The Hilbert basis of $\overline{S_A}$ is given by $\mathcal{H}(\overline{S_A})=\{h'_1,\dots,h'_7\}$, where $h'_i$ corresponds to the $i$-th column of the following matrix
$$\left[\begin{array}{rrrrrrr}
0 & 0 & 0 & 1 & 1 & 0 & 1 \\
0 & 1 & 2 & 0 & 1 & 1 & 1 \\
0 & 0 & 3 & 0 & 3 & 1 & 2 \\
1 & 0 & 0 & -1 & 0 & 0 & 0
\end{array}\right].
$$
Notice that $\operatorname{det}_3(h'_1\,h'_2\,h'_4\,h'_6)=1$. Let $A'=\{h'_1,h'_2,h'_4,h'_6\}$. The semigroup $S_{A'}$ generated by $\mathcal{G}_{A'}$ is pointed and so  $X_2 = X(\overline{S_{A'}})$ is an affine chart of the normalized Nash blowup of $X_1$. We denote by $\omega_{3,2}$ the cone generated by $\overline{S_{A'}}$ in $M_\R$. Then, the cone $\omega_{3,2}$ is generated by the columns of the matrix
$$B_{3,2}=\left[\begin{array}{rrrrr}
0 & 0 & 0 & 1 & 1 \\
0 & 1 & 1 & 0 & 0 \\
0 & 0 & 3 & 0 & 3 \\
1 & 0 & 0 & -1 & 0
\end{array}\right]\,.$$
We conclude that $X_2$ is an affine chart of the second iteration of the normalized Nash blowup of $X$. Now we have that $X$ is isomorphic to $X_2$. Indeed, the automorphism $U\colon M\to M$ given by the unimodular matrix
$$U=\left[\begin{array}{rrrr}
1 & 1 & -1 & 0 \\
0 & 0 & 0 & 1 \\
3 & 0 & -1 & 3 \\
0 & -1 & 1 & 0
\end{array}\right]\,.$$
provides an isomorphism from $\omega_{3,0}$ to $\omega_{3,2}$, concluding the proof of the theorem. \qed

\begin{remark} \label{rem:non-alg-closed}
It is natural to ask whether the previous counterexamples also hold over non-algebraically closed fields. The first difficulty towards that question relies on the fact that we have a general combinatorial description for the Nash blowup of toric varieties only over algebraically closed fields. Indeed, \cite[Proposition 60]{GoTe14} and \cite[Theorem 1.9]{DJNB} require that hypothesis. Hence, we cannot apply the combinatorial description we used in the previous proofs. However, at least over characteristic zero fields, there is another approach that allows us to show that the counterexample also holds in that case. The following discussion is based on \cite[Construction 4.4]{GrMi12}. We apply it to the first example of \cref{sec:proof-counter} (recall the notation of \cref{lem S char 0}).

Let $X(S)\subset\K^7$ be the corresponding 4-dimensional toric variety. Recall that $X(S)$ can also be obtained as the Zariski closure of the image of the monomial map
$$\phi:(\K^*)^4\to\K^7,\,\,\,x\mapsto (x^{h_1},\ldots,x^{h_7}).$$
The tangent space $T_{\phi(x)}X(S)$ is determined by the Jacobian matrix  $J=J(x^{h_1},\ldots,x^{h_7})$. Using the Pl\"ucker embedding $P:\operatorname{Grass}(4,7)\hookrightarrow\mathbb{P}_{\K}^{\binom{7}{4}-1}$, the Nash blowup of $X(S)$ is then obtained as the Zariski closure of 
$$\{(\phi(x),(\ldots:\Delta_{i_1i_2i_3i_4}:\ldots))\mid x\in(\K^*)^4\}\subset X(S)\times\mathbb{P}_{\K}^{\binom{7}{4}-1},$$
where $\Delta_{i_1i_2i_3i_4}$ denotes the maximal minor defined by the rows $i_1,i_2,i_3,i_4$ of $J$.

Consider the coordinate of $\mathbb{P}_{\K}^{\binom{7}{4}-1}$ corresponding to the position of $\Delta_{1,2,3,5}$. Dividing the other coordinates by this minor we obtain an affine chart of $X(S)^*$. This affine chart turns out to be isomorphic to $X(S)$, just like in the proof in \cref{sec:proof-counter}. This can be verified through straightforward computations: it is just a matter of computing the 35 minors of the Jacobian matrix (all of them monomials) and dividing them by $\Delta_{1,2,3,5}$; these quotients being monomials, their exponents generate a semigroup which, in turn, is generated by the same set $H$ from \cref{sec:proof-counter}. Notice that it is in the division by $\Delta_{1,2,3,5}$ that the characteristic zero assumption is used. 
\end{remark}

\section{Remarks on computer experimentation} \label{sec:computer}

Even though the theorem was proven without the need for computer assistance, the counterexamples presented in the proof were obtained through computer experimentation. We now describe the process that led to their discovery. A complete exposition of the implementation, optimizations, and computational results is more naturally suited to an independent presentation and will be provided in a forthcoming paper~\cite{CDLLcomputational}.

We implemented the algorithm for the normalized Nash blowup using SageMath \cite{sagemath}, which includes built-in methods for handling cones, fans, and normal toric varieties. Let $\omega \subset M_\mathbb{R}$ be a pointed rational polyhedral cone, and let $S = \omega \cap M$. The normalized Nash blowup of the normal affine toric variety $X(S)$ is computed as follows:
\begin{enumerate}
    \item Compute the Hilbert basis $\mathcal{H}$ of the semigroup $S$.
    \item Construct the polyhedron $\mathcal{N}_p(S)\subset M_\mathbb{R}$ as defined in~\eqref{eq:norm-Nash}.
\end{enumerate}
Then, as stated below  \eqref{eq:norm-Nash}, a  set of covering affine charts of the normalized Nash blowup of $X(S)$ is given by $X(S_v)$, where $v$ runs over all the vertices of the polyhedron  $\mathcal{N}_p(S)$, $\omega_v = \operatorname{Cone}\left(\mathcal{N}_p(S)-v\right)$, and $S_v=\omega_v\cap M$.

Implementing this algorithm in SageMath is straightforward. For the first step, SageMath provides a built-in function to compute the Hilbert basis of the semigroup associated with a cone. The second step only involves computing determinants in low dimensions and determining the vertices of polyhedra, functionalities that are also readily available in SageMath. We then iterate this algorithm, discarding affine charts that are already smooth. The affine toric variety $X(S)$ is resolved by successive normalized Nash blowups if and only if this iteration eventually terminates.

\medskip

Our initial motivation for implementing this algorithm was to investigate termination in the three-dimensional case and characteristic 0. Let $\rho(d,n)$ be the cone in $M_\mathbb{R}=\mathbb{R}^d$ generated by the columns of the $d \times d$ matrix
\begin{align*}
\left[
\begin{array}{c|c}
I_{d-1} & \mathbf{1} \\
\hline
\mathbf{0} & n
\end{array}
\right]
\end{align*}
We began investigating the family of normal toric varieties $X(S)$ whose semigroups are given by $S = \rho(3,n) \cap M$ for increasing values of $n$. 
These cones were not chosen at random; in fact, they are the source of pathological examples in the theory of lattices polytopes; see for instance \cite{reeve1957volume} and \cite{beck2015very}. 
Also, they are a rare example of a non-unimodular cone where we can nevertheless write down the whole Hilbert basis.
At this stage, we observed that these toric varieties yield particularly interesting normalized Nash resolutions. 
Although they require several steps to resolve, no clear pattern emerges in the behavior of the resolution process. 
We computed the full normalized Nash resolution for the cones $\rho(3,n)$ with $2 \leq n \leq 1000$, all of which lead to termination.

After repeated unsuccessful attempts to uncover patterns in the three-dimensional case that could shed light on the appropriate invariants to prove the termination of the algorithm, we eventually decided to shift our focus. 
We began drafting our paper~\cite{CDLLcharfree}, which presents the positive results we were able to obtain. 
Our initial computer implementation of the normalized Nash blowup was restricted to the three-dimensional setting. 
However, while preparing the examples for that paper, we decided to generalize the implementation to arbitrary dimension.
Naturally, we then used this extended version to investigate the four-dimensional case.

\medskip

Analogous to the three-dimensional case, we then studied the cones $\rho(4,n)$ for increasing values of~$n$. 
It is in this context that counterexamples arose. 
Indeed, the normal toric varieties $X(S)$ whose semigroups are given by $S = \rho(4,n) \cap M$ for $n = 2, 3, 4$ all lead to the termination of the algorithm. 
However, for  the cone $\rho(4,5)$, generated by the columns of the matrix  
$$\left[\begin{array}{rrrr}
1 & 0 & 0 & 1 \\
0 & 1 & 0 & 1 \\
0 & 0 & 1 & 1 \\
0 & 0 & 0 & 5
\end{array}\right]\, ,$$
the normal toric variety $X(S)$ fails to resolve via normalized Nash blowup. Crucially, the number of non-smooth cones stabilizes at 744 after the 13th iteration. It was examining these 744 singular charts, of which there were only 34 pairwise non-isomorphic, that we found the cone determined by the columns of $B$ in \eqref{eq:count-char-0}, giving place to the counterexample in characteristics different from 2 or 3. Moreover, further investigation revealed that the chart appears for the first time in the 7th iteration. Regarding characteristics 2 and 3 we had to search for a different counterexample. Remarkably, the iteration of the normalized Nash blowup of the toric variety $X$ corresponding to the cone $\rho(4,5)$ also fails to resolve in characteristics 2 and 3.

In characteristic $2$, an isomorphic copy of the toric variety $X(S)$, where $S = \omega_2 \cap M$ and $\omega_2$ is the cone generated by the columns of $B_2$ in~\eqref{eq:count-char-2}, appears for the first time in the 3rd iteration of the normalized Nash blowup of $X$. The toric variety $X(S)$ is resolved by iteration of the normalized Nash blowup in all characteristics different from $2$.

In characteristic $3$, an isomorphic copy of the toric variety $X(S)$, where $S = \omega_{3,0} \cap M$ and $\omega_{3,0}$ is the cone generated by the columns of $B_{3,0}$ in~\eqref{eq:count-char-3}, appears for the first time in the 4th iteration of the normalized Nash blowup of $X$. The toric variety $X(S)$ is resolved by iteration of the normalized Nash blowup in all characteristics different from $2$ and $3$. In the case of characteristic $2$, the second iteration of the normalized Nash blowup contains an affine chart isomorphic to the counterexample presented in the proof. 

\medskip

After extensive computer experimentation, we were unable to identify other four-dimensional toric varieties exhibiting the property of being isomorphic to an affine chart appearing in an iteration of their own normalized Nash blowup. 
All counterexamples we have found ultimately fall into the loops described by the toric varieties in~\cref{sec:proof-counter}. 
This suggests that the counterexamples presented in this article are extremely exceptional and merit further study from multiple perspectives.

\medskip

Finally, we turned our attention to the five-dimensional case, which is the highest dimension currently feasible in our implementation on standard desktop computers. 
This case behaves quite differently. 
In characteristic zero, we have been able to find several independent counterexamples. So far, we have identified 11 disjoint loops. 
Specifically, we have found singular normal toric varieties in which an affine chart isomorphic to the original variety appears for the first time in the $n$-th iteration of the normalized Nash blowup, for $n = 1, 2, 4, 5, 8,$ and $9$.

\bibliographystyle{alpha}
\bibliography{ref}

\end{document}